\documentclass[11pt]{amsart}
\usepackage{amsmath,amssymb,latexsym,amsfonts,amscd}
\newtheorem{theorem}{Theorem}[section]
\newtheorem{proposition}[theorem]{Proposition}
\newtheorem{definition}[theorem]{Definition}

\newtheorem{corollary}[theorem]{Corollary}

\begin{document}

\title[$\mathcal D$-modules, characteristic-free]{A characteristic-free proof of a basic result on $\mathcal D$-modules}
\author{Gennady Lyubeznik}
\address{Dept. of Mathematics, University of Minnesota, Minneapolis,
MN 55455}
\email{gennady@math.umn.edu}
\thanks{NSF support through grants DMS-0202176 and DMS-0701127 is gratefully acknowledged. }

\begin{abstract} Let $k$ be a field, let $R$ be a ring of polynomials in a finite number of variables over $k$, let $\mathcal D$ be the ring of $k$-linear differential operators of $R$ and let $f\in R$ be a non-zero element. It is well-known that $R_f$, with its natural $\mathcal D$-module structure,  has finite length in the category of $\mathcal D$-modules. We give a characteristic-free proof of this fact. To the best of our knowledge this is the first characteristic-free proof.  \end{abstract}

\maketitle

\section{Introduction.} Throughout this paper $k$ is a field, $R=k[x_1,\dots,x_n]$ is the ring of polynomials in a finite number of variables over $k$ and $\mathcal D$ is the ring of $k$-linear differential operators of $R$. The natural $\mathcal D$-action on $R$ induces a $\mathcal D$-module structure on $R_f$ for every $0\ne f\in R$. The goal of this paper is to give a characteristic-free proof of the following well-known fact.
\begin{theorem}\label{main} 
$R_f$ has finite length in the category of $\mathcal D$-modules.
\end{theorem}
In characteristic 0 this is due to J. Bernstein \cite{[Be1], [Be2]} and in characteristic $p>0$ to R. B\o gvad \cite{[Bo]}. In both cases proofs are based on suitable notions of holonomicity but the definitions of holonomicity in each of these two cases are completely different. 

Our characteristic-free proof is made possible by V. Bavula's wonderful paper \cite{Ba} where a characteristic-free definition of holonomic modules is given. But the focus of \cite{Ba} is the characteristic $p>0$ case and this assumption is routinely made in the statements and used in the proofs. 

In this paper we simplify and characteristic-freeify those of Bavula's results that are needed for a proof of Theorem \ref{main}.

Finiteness properties of local cohomology modules for regular rings containing a field had originally been proven by two completely different methods in characteristic $p>0$ \cite{[HuSh]} and in characteristic 0 \cite{[Ly1]}. In \cite{FinProp} we used $\mathcal D$-modules to give proofs of these finiteness properties that are characteristic-free modulo the fact that $\mathcal R_f$, where $\mathcal R=k[[x_1,\dots,x_n]]$ is the ring of formal power series in a finite number of variables over $k$, has finite length in the category of $k$-linear $\mathcal D$-modules of $\mathcal R$. The proofs of this complete local analogue of Theorem \ref{main} are still completely different in characteristic 0 \cite{Bjork} and in characteristic $p>0$ \cite{F-mod}. 

Our proof of Theorem \ref{main} leads to a characteristic-free proof of the finiteness properties of local cohomology modules over polynomial rings. And it suggests a way to find a similar proof in general, i.e. for all regular local rings containing a field: through a suitable characteristic-free definition of holonomicity in the complete local case that would lead to a proof of an analogue of Theorem \ref{main} in this case. Such a definition is yet to be discovered.

This paper is self-contained.

\section{Preliminaries.} Let $D_{t,i}=\frac{1}{t!}\frac{\partial^t}{\partial  x^t_i}:R\to R$ be the $k[x_1\dots,x_{i-1},x_{i+1},\dots,x_n]$-linear map that sends $x_i^v$ to $\binom{v}{t}x^{t-v}$ ($D_{0,i}$ is the identity map). Even though $\frac{1}{t!}$ is part of the notation, $D_{t,i}$ exists in all characterisitcs because $\binom{v}{t}$ is an integer. 

The ring $R$ is in a natural way a subring of End$_kR$ (every element of $R$ corresponds to the multiplication by that element on $R$) and the following equality holds in End$_kR$.

\begin{proposition}\label{formula}
$D_{t,i}\cdot f=\Sigma_{s=0}^{s=t}D_{s,i}(f)\cdot D_{t-s,i}$ for every $f\in R$.
\end{proposition}

\emph{Proof.} We have to show that for every $g\in R$
$$D_{t,i}(f\cdot g)=\Sigma_{s=0}^{s=t}D_{s,i}(f)\cdot D_{t-s,i}(g)$$
which is the well-known formula for the higher derivative of a product $$\frac{\partial^t}{\partial  x^t_i}(f\cdot g)=\Sigma_{s=0}^{s=n}\binom{t}{s}\frac{\partial^s}{\partial  x^s_i}(f)\cdot \frac{\partial^{t-s}}{\partial  x^{t-s}_i}(g)$$ divided by $t!$ \qed

\begin{corollary}\label{relat}
(a) $D_{t,i}$ commutes with $x_j$ for $j\ne i$ and with all $D_{s,j}$.

(b) $D_{t,i}x^w_i=\Sigma_{s=0}^{s=t}\binom{w}{s}x^\ell_iD_{t-s,i}$.

(c) $D_{t,i}\cdot D_{s,j}=\binom{s+t}{s}D_{t+s,i}$.
\end{corollary}

\emph{Proof.} (a) and (c) are straightforward and (b) is \ref{formula} with $f=x_i^w$.\qed

\medskip

The ring $\mathcal D$ of $k$-linear differential operators of $R$ is the $k$-subalgebra of End$_kR$ generated by $R$ and all the $D_{t,i}$s. Corollary \ref{relat} implies that 
the products $\{x_1^{i_1}\cdots x_n^{i_n}\cdot D_{t_1,1}\cdots D_{t_n,n}\}$ where $i_1,\dots, i_n,t_1,\dots,t_n$ range over all the $2n$-tuples of non-negative integers, are a $k$-basis of $\mathcal D$. Indeed, every element of $\mathcal D$ is by definition a linear combination of products of $D_{t,i}$s and $x_j$s. Using relations \ref{relat}(a)-(c) we can write every such product as a linear combination of products of the form $x_1^{i_1}\cdots x_n^{i_n}\cdot D_{t_1,1}\cdots D_{t_n,n}$. Thus $\mathcal D$ is free left $R$-module on the products $D_{t_1,1}\cdots D_{t_n,n}$ and similarly, it is a free right $R$-module on these same products.

\begin{corollary}\label{socle}
$x^w\cdot D_{t,i}\in \mathcal Dx_i$ if $w>t$ and $x^t\cdot D_{t,i}-(-1)^t\in \mathcal Dx_i$.
\end{corollary}

\emph{Proof}. Isolating $x_i^w\cdot D_{t,i}$ from \ref{relat}b we get
$$x_i^w\cdot D_{t,i}=D_{t,i}\cdot x_i^w-\Sigma_{s=1}^{s=t}x_i^{w-s}\cdot D_{t-s,i}$$ which implies both containments by induction on $t$. \qed

\begin{proposition}\label{D/Dm}
Let $\mathfrak m\subset R$ be a $k$-rational maximal ideal of $R$ (this means that the natural map $k\hookrightarrow R/\mathfrak m$ is bijective). If $\delta\in \mathcal D$, we denote $\bar\delta\in \mathcal D/\mathcal D\mathfrak m$ the image of $\delta$ under the natural surjection $\mathcal D\to \mathcal D/\mathcal D\mathfrak m$.

(i) $\mathcal D/\mathcal D\mathfrak m$ is the $k$-vector space with basis $\{\overline{D_{t_1,1}\cdots D_{t_n,n}}\}$ as $t_1,\dots,t_n$ range over all non-negative integers.

(ii)  Every element of $\mathcal D/\mathcal D\mathfrak m$ is annihilated by a power of $\mathfrak m$ and the socle of $\mathcal D/\mathcal D\mathfrak m$ is generated by $\bar 1$.
\end{proposition}

\emph{Proof.} (i) follows from the fact that $\mathcal D$ is a free right $R$-module on the products $D_{t_1,1}\cdots D_{t_n,n}$ and $\mathcal D/\mathcal D\mathfrak m=\mathcal D\otimes_R(R/\mathfrak m)$.

(ii) Since the natural map $k\hookrightarrow R/\mathfrak m$ is bijective, $\mathfrak m=(x_1-c_1,\dots,x_n-c_n)$ where $c_1,\dots,c_n\in K$. Viewing $x_j-c_j$ as a new $x_j$ we can assume that $\mathfrak m=(x_1,\dots,x_n)$. According to \ref{socle}, $x_i^{t_i+1}$ annihilates $\{\overline{D_{t_1,1}\cdots D_{t_n,n}}\}$, hence every element of $\mathcal D/\mathcal D\mathfrak m$ is annihilated by a power of $\mathfrak m$. Clearly $\bar 1$ belongs to the socle. It remains to show that every non-zero element $z$ can be sent to the socle by multiplication by an element of $R$. According to (i) $z$ is a $k$-linear combination of a finite number of $\overline{D_{t_1,1}\cdots D_{t_n,n}}$. Let $D_{t_1,1}\cdots D_{t_n,n}$ have a maximal $t_1+\dots +t_n$ among all the $\overline{D_{t_1,1}\cdots D_{t_n,n}}$ with non-zero coefficients in this linear combination. Hence for every other $\overline{D_{t'_1,1}\cdots D_{t'_n,n}}$ with non-zero coefficient in the linear combination there is $j$ such that $t_j>t'_j$. It follows from \ref{relat} that $x_j^{t_j}D_{t',j}\in \mathcal D\mathfrak m$. Hence $x_1^{t_1}\cdots x_n^{t_n}\overline{D_{t'_1,1}\cdots D_{t'_n,n}}=0$. It similarly follows from \ref{relat} that $x_1^{t_1}\cdots x_n^{t_n}\overline{D_{t'_1,1}\cdots D_{t'_n,n}}=(-1)^{t_1+\dots+t_n}\bar 1$. Hence $(-1)^{t_1+\dots+t_n}x_1^{t_1}\cdots x_n^{t_n}z=\bar 1$. \qed

\begin{corollary}\label{lin.ind}
Let $\mathfrak m\subset R$ be a maximal ideal such that $R/\mathfrak m$ is a finite separable field extension of $k$. Let $M$ be a $\mathcal D$-module and let $z\in M$ be a non-zero element such that  its annihilator in $R$ is $\mathfrak m$. The set $\{D_{t_1,1}\cdots D_{t_n,n}z\}$, as $t_1,\dots,t_n$ range over all non-negative integers, is linearly independent over $k$.
\end{corollary}

\emph{Proof.} Replacing $M$ by the $\mathcal D$-submodule generated by $z$ we can assume that $M$ is generated by $z$. Let $K$ denote the algebraic closure of $k$, let $R'=K\otimes_kR=K[x_1,\dots,x_n]$, $\mathcal D'=K\otimes_k\mathcal D$ and $M'=K\otimes_kM$. Then $\mathcal D'$ is the ring of $K$-linear differential operators of $R'$ and $M'$ is naturally a $\mathcal D'$-module.  Identifying $M$ with the subset $1\otimes_kM$ of $M'$ we conclude that it is enough to show that the set $\{D_{t_1,1}\cdots D_{t_n,n}z\}\subset M'$ is linearly independent over $K$.

Let $\mathfrak m_1\dots,\mathfrak m_s$ be the maximal ideals of $R'$ that lie over $\mathfrak m$. Since the field extension $k\hookrightarrow R/\mathfrak m$ is separable, $K\otimes_kR/\mathfrak m$ is reduced. Therefore $K\otimes_kR/\mathfrak m=R'/(\cap_i\mathfrak m_i)$. This implies that $R'z=K\otimes_kRz\cong R'/(\cap_i\mathfrak m_i)$ since $Rz\cong R/\mathfrak m$. Now $M'$ being generated by $z$ is a surjective image of $\mathcal D'/\mathcal D'(\cap_i\mathfrak m_i)$  via the surjection is $\mathcal D'/\mathcal D'(\cap_i\mathfrak m_i)\stackrel{\bar 1\mapsto z}{\to}M'$. 

But $\mathcal D'/\mathcal D'(\cap_i\mathfrak m_i)=\mathcal D'\otimes_{R'}R'/(\cap_i\mathfrak m_i)=\mathcal D'\otimes_{R'}(\oplus_iR'/\mathfrak m_i)=\oplus_iD'/D'\mathfrak m _i.$ According to \ref{D/Dm}, the socle of each $D'/D'\mathfrak m _i$ is generated by $\bar 1$, hence so is the socle of $\mathcal D'/\mathcal D'(\cap_i\mathfrak m_i)$. This means the surjection induces a bijection on the socles and therefore it is itself a bijection. Thus $\mathcal D'/\mathcal D'(\cap_i\mathfrak m_i)$ is isomorphic to $M$ via an isomorphism that sends $\overline{D_{t_1,1}\cdots D_{t_n,n}}$ to $\{D_{t_1,1}\cdots D_{t_n,n}z\}$. But the set $\{\overline{D_{t_1,1}\cdots D_{t_n,n}}\}$ is linearly independent (this is a consequence of \ref{D/Dm} after a natural projection onto some $D'/D'\mathfrak m _i$). \qed

\begin{definition}
The Bernstein filtration $k=\mathcal F_0\subset \mathcal F_1\subset \mathcal F_2\subset\dots$ on $\mathcal D$ is defined by setting $\mathcal F_s$ to be the $k$-linear span of the set of products $\{x_1^{i_1}\cdots x_n^{i_n}\cdot D_{t_1,1}\cdots D_{t_n,n}|i_1+\dots i_n+t_1+\dots+t_n\leq s\}$. \end{definition}

It follows from \ref{relat} that $\mathcal F_i\cdot\mathcal F_j\subset \mathcal F_{i+j}$.

\section{Proof of Theorem \ref{main}} 

The technical heart of our proof is the following proposition.

\begin{proposition}\label{geq}
(cf. \cite[9.3]{Ba}) Assume the field $k$ is separable. Let $M$ be a $\mathcal D$-module and let $z\in M$ be an element such that the annihilator of $z$ in $R$ is a prime ideal of $R$. Then ${\rm dim}_k(\mathcal F_iz)\geq \binom{n+i}{i}$.
\end{proposition}

\emph{Proof.} Let $d={\rm dim}R/P$, let $h=n-d$, and let $K$ be the fraction field of $R/P$. Since the transcendence degree of $K$ over $k$ equals $d$ and $k$ is separable, after a possible permutation of indices we can assume that $x_{h+1},\dots,x_n$ are algebraically independent over $k$ in $K$ and $K$ is finite and separable over the field of rational functions $\mathcal K=k(x_{h+1},\dots,x_{n})$. 

Let $R'=\mathcal K\otimes_RR=\mathcal K[x_1,\dots,x_h]$, let $\mathcal D'$ be the ring of $\mathcal K$-linear differential operators of $R'$ and let $M'=\mathcal K\otimes_RM$. The ring $\mathcal D'$ is a free $R'$-module on the products $D_{t_1,1}\cdots D_{t_h,h}$. Since each such product commutes with $x_j$ for $j>h$, its action on $M$ naturally extends to an action on $M'$ making $M'$ a $\mathcal D'$-module. It follows from \ref{lin.ind} that the set of elements $\{D_{t_1,1}\cdots D_{t_h,h}z\}\subset M'$, as $t_1,\dots,t_h$ run through all non-negative integers, is linearly independent over $\mathcal K$. Setting $R''=k[x_{h+1},\dots,x_{n}]$ ($\mathcal K$ is the fraction field of $R''$) we conclude that the sum $\Sigma_{t_1,\dots,t_{h}} R''D_{t_1,1}\cdots D_{t_h,h}z$ of $R''$-submodules of $M$ is direct, i.e. the natural surjective $R''$-module map $\oplus_{t_1,\dots,t_{h}} R'D_{t_1,1}\cdots D_{t_h,h}\to \Sigma_{t_1,\dots,t_{h}} R''D_{t_1,1}\cdots D_{t_h,h}z$ that sends every product $D_{t_1,1}\cdots D_{t_h,h}\in \oplus_{t_1,\dots,t_{h}} R''D_{t_1,1}\cdots D_{t_h,h}$ to $D_{t_1,1}\cdots D_{t_h,h}z\in M$ is bijective. And this implies that the set $\{x_{h+1}^{i_{h+1}}\cdots x_n^{i_n}D_{t_1,1}\cdots D_{t_h,h}z\}$ of elements of $M$, as $t_1,\dots,t_h,i_{h+1},\dots,i_n$ run over all non-negative integers, is linearly independent over $k$.

The elements of this set with $t_1+\cdots t_h+i_1+\cdots+i_n\leq i$ belong to $\mathcal F_iz$. The number of these elements equals the number of monomials of degree at most $i$ in $n$ variables which is well-known to equal $\binom{n+i}{i}$.\qed

\begin{definition}
A $k$-filtration on a $\mathcal D$-module $M$ is an ascending chain of $k$-vector spaces $M_0\subset M_1\subset M_2\subset\dots$ such that $\cup_iM_i=M$ and $\mathcal F_iM_j\subset M_{i+j}$ for all $i$ and $j$.
\end{definition}

For example, the Bernstein filtration on $\mathcal D$ is a $k$-filtration.

\begin{corollary}
Let $M_0\subset M_1\subset M_2\subset\dots$ be a $k$-filtration on a $\mathcal D$-module $M$. There exists an integer $j$ such that ${\rm dim}_kM_i\geq \binom{n+i-j}{i-j}$ for all $i\geq j$.
\end{corollary}

\emph{Proof.} Let $k'$ be the algebraic closure of $k$ and $R'=k'\otimes_kR$. Then the ring of $k'$-linear differential operators of $R'$ is $\mathcal D'=k'\otimes_k\mathcal D$ and $M'=k'\otimes_RM$ is in a natural way a $\mathcal D'$-module with $k'$-filtration $M'_0\subset M'_1\subset M'_2\subset\dots$ where $M_i'=k'\otimes_kM_i$. Since ${\rm dim}_kM_i={\rm dim}_{k'}M'_i$, we can and do assume that $k$ is algebraically closed and in particular separable.

Let $P\subset R$ be an associated prime ideal of $M$ in $R$. This means there exists an element $z\in M$ such that the annihilator of $z$ in $R$ is $P$. Let $j$ be the smallest integer such that $z\in M_j$. Clearly $M_i\supset \mathcal F_{i-j}z$, so we are done by \ref{geq}\qed

\medskip

The following definition of holonomicity is equivalent to but somewhat simpler than Bavula's original definition \cite[pp. 185, 198]{Ba}; in particular we do not require the module $M$ to be finitely generated. But Theorem \ref{FinLen} implies that every holonomic module is finitely generated (this fact is not used in the sequel). 

\begin{definition}
A $\mathcal D$-module $M$ is holonomic if it has a $k$-filtration with ${\rm dim}_kM_i\leq Ci^n$ for all $i$ where $C$ is a constant independent of $i$.
\end{definition}

\begin{theorem}\label{FinLen}
(cf.  \cite[9.6]{Ba}) Every holonomic $\mathcal D$-module has finite length in the category of $\mathcal D$-modules. In fact if $M_0\subset M_1\subset\dots$ is a $k$-filtration on $M$ with ${\rm dim}_kM_i\leq Ci^n$, then the length of $M$ in the category of $\mathcal D$-modules is at most $n!C$.
\end{theorem}

\emph{Proof.} Let $0=M^0\subset M^1\subset \dots M^{\ell}=M$ be a filtration of $M$ in the category of $\mathcal D$-modules. Then $(M^s/M^{s-1})_i=(M_i\cap M^s)/(M_i\cap M^{s-1})$ is a $k$-filtration on the $\mathcal D$-module $M^s/M^{s-1}$. Hence there is an integer $j_s$ such that dim$_k(M^s/M^{s-1})_i\geq\binom{n+i-j_s}{i-j_s}$ for all $i\geq j_s$. 

But $M_i=\oplus_{s=1}^{s=\ell}(M^s/M^{s-1})_i$ because these are vector spaces over a field $k$. Therefore dim$_kM_i=\Sigma_{s=1}^{s=\ell}{\rm dim}_k(M^s/M^{s-1})_i\geq \Sigma_{s=1}^{s=\ell}\binom{n+i-j_s}{i-j_s}$ for all sufficiently big $i$. This implies $Ci^n\geq \Sigma_{s=1}^{s=\ell}\binom{n+i-j_s}{i-j_s}$ for all sufficiently big $i$.

But $\binom{n+i-j_s}{i-j_s}$, for fixed $n$ and $j_s$, is a polynomial in $i$ of degree $n$ and top coefficient $\frac{1}{n!}$. Hence $p(i)=\Sigma_{s=1}^{s=\ell}\binom{n+i-j_s}{i-j_s}$ is a polynomial in $i$ of degree $n$ and top coefficient $\frac{\ell}{n!}$. The inequality $Ci^n\geq p(i)$ holds for all sufficiently big $i$ only if $C\geq \frac{\ell}{n!}$, i.e. $\ell\leq n!C$. \qed

\medskip

If $M$ is a $\mathcal D$-module and $f\in R$ is a non-zero element, the module $M_f$ acquires a structure of $\mathcal D$-module as follows. The formula \ref{formula} implies 
$$f\cdot D_{t.i}=D_{t,i}\cdot f-\Sigma_{s=1}^{s=n}D_{s,i}(f)\cdot D_{t-s,i}.$$ Replacing $f$ by $f^j$ in this equality and then applying it to $\frac{m}{f^j}\in M_f$ and multiplying on the left by $f^{-j}$ we get 

\begin{gather}
D_{t.i}(\frac{m}{f^j})=f^{-j}\cdot D_{t,i}(m)-\Sigma_{s=1}^{s=n}f^{-j}\cdot D_{s,i}(f^j)\cdot D_{t-s,i}(\frac{m}{f^j}) 
\end{gather}

This leads to a definition of the action of $D_{t,i}$ on $M_f$ by induction on $t$ the case $t=0$ being trivial (since $D_{0,i}$ is the identity map).

Modules of type $M_f$ are not considered in \cite{Ba}.

\begin{corollary}\label{fin.length}
If $M$ is a holonomic module and $f\in R$, then $M_f$ is holonomic.
\end{corollary}

\emph{Proof.} Let $M_0\subset M_1\subset \dots$ be a $k$-filtration of $M$ with dim$_kM_i\leq Ci^n$. Let $d$ be the degree of $f$ and let $M'_0\subset M'_1\subset\dots$ be the filtration on $M_f$ defined by $M'_i=\{\frac{m}{f^i}|m\in M_{i(d+1)}\}$. We claim this is a $k$-filtration, i.e. $\cup_iM'_i=M_f$ and $\mathcal F_iM'_j\subset M'_{i+j}$ for all $i$ and $j$.

Indeed, let $\frac{m}{f^w}\in M_f$ be any element. Assume $m\in M_u$. If $u\leq w(d+1)$, then $m\in M_{w(d+1)}$ hence $\frac{m}{f^w}\in M'_w.$ If $u> w(d+1)$ let $v=u-w(d+1)$. Since $f^v\in \mathcal F_{vd}$, it follows that $f^v\cdot m\in M_{vd+u}$. Since $vd+u=(v+w)(d+1)$, we conclude that  $\frac{m}{f^w}=\frac{f^v\cdot m}{f^{w+v}}\in M'_{v+w}$.  This shows that $\cup_iM'_i=M_f$.

To prove that $\mathcal F_iM'_j\subset M'_{i+j}$ all we need to show is that $x_uM'_j\subset M'_{j+1}$ and $D_{t, u}M'_j\subset M'_{t+j}$ for every $u\in \{1,2,\dots, n\}$. Let $\frac{m}{f^j}\in M'_j$ where $m\in M_{j(d+1)}$. 

Since $x_u\cdot f\in \mathcal F_{d+1}$, it follows that $(x_u\cdot f)m\in M'_{(j+1)(d+1)}$. Therefore $x_u\cdot \frac{m}{f^j}=\frac{(x_u\cdot f)m}{f^{j+1}}\in M'_{j+1}$. This shows that $x_uM'_j\subset M'_{j+1}$.

To prove that $D_{t, u}M'_j\subset M'_{t+j}$ we use induction on $t$ the case $t=0$ being trivial since $D_{0,u}$ is the identity map. It is enough to show that all the terms on the right side of (1), i.e. $f^{-j}\cdot D_{t,u}(m)$ and $f^{-j}\cdot D_{s,u}(f^j)\cdot D_{t-s,u}(\frac{m}{f^j})$, where $s\geq 1$, belong to $M'_{t+j}$.

Since $f\in \mathcal F_d$, $D_{t,u}\in \mathcal F_t$ and $m\in M_{j(d+1)}$, it follows that $f^t\cdot D_{t,u}\in \mathcal F_{td+t}$ and $f^t\cdot D_{t,u}(m)\in M_{j(d+1)+td+t=(t+j)(d+1)}$. Thus $f^{-j}\cdot D_{t,u}(m)=\frac{f^t\cdot D_{t,u}(m)}{f^{t+j}}$ belongs to $M'_{t+j}$ because the top of this fraction belongs to $M_{(t+j)(d+1)}$. 

If $s\geq 1$, then $D_{t-s,u}(\frac{m}{f^j})\in M'_{t-s+j}$ by the induction hypothesis, i.e. there exists $m_{t-s}\in M_{(t-s+j)(d+1)}$ such that $D_{t-s,u}(\frac{m}{f^j})=\frac{m_{t-s}}{f^{t-s+j}}$. It follows by induction on $j$ using formula \ref{formula} that $D_{s,i}(f^j)$ is divisible by $f^{j-s}$, i.e. $D_{s,i}(f^j)=f^{j-s}\cdot f_s$. Hence $f^{-j}\cdot D_{s,u}(f^j)\cdot D_{t-s,u}(\frac{m}{f^j})=\frac{f_s\cdot m_{t-s}}{f^{t+j}}$. Since the polynomial $D_{s,i}(f^j)$ has degree $dj-s$, the polynomial $f_s$ has degree $ds-s$. Hence $f_s\cdot m_{t-s}\in M_{ds-s+(t-s+j)(d+1)}\subset M_{(t+j)(d+1)}$. The latter containment is because $ds-s+(t-s+j)(d+1)\leq (t+j)(d+1)$. This shows that $f^{-j}\cdot D_{s,u}(f^j)\cdot D_{t-s,u}(\frac{m}{f^j})\in M'_{t+j}$ and completes the proof that $D_{t, u}M'_j\subset M'_{t+j}$ which in turn completes the proof of our claim that $M'_0\subset M'_1\subset\dots$ is a $k$-filtration. 

Clearly, dim$_kM'_i\leq{\rm dim}_kM_{i(d+1)}\leq C(i(d+1))^n$. This implies that dim$_kM'_i\leq C'i^n$ where $C'=C(d+1)^n$.\qed

\medskip

The filtration by degree on the $\mathcal D$-module $M=R$ (i.e. $M_i$ consists of all the polynomials of degree at most $i$) shows that $R$ is a holonomic $\mathcal D$-module. Now Theorem \ref{main} follows from \ref{fin.length}. \qed

\medskip

\section{Some open problems}

1. Let $0\to M'\to M\to M''\to 0$ be an exact sequence in the category of $\mathcal D$-modules. It is not hard to see that if $M$ is holonomic, then so are $M'$ and $M''$. In characteristic 0 the converse is also true, i.e. if $M'$ and $M''$ are holonomic, then so is $M$. Is this true in characteristic $p>0$ as well?

\medskip

2. Let $M$ be a holonomic $\mathcal D$-module. Since $M$ has finite length, it is finitely generated as a $\mathcal D$-module. This implies that there is a $k$-filtration $M_0\subset M_1\subset \dots$ on $M$ such that $M_0$ is finite-dimensional over $k$ and $M_i=\mathcal F_iM_0$ (just take $M_0$ to be the $k$-span of a finite set of $\mathcal D$-generators of $M$). It is not hard to show that ${\rm lim sup}_{n\to \infty}\frac{{\rm dim}_kM_i}{i^n}$ is independent of $M_0$. It is well-known that lim$_{n\to \infty}\frac{{\rm dim}_kM_i}{i^n}$ exists in characteristic 0 and, moreover, $n!({\rm lim}_{n\to \infty}\frac{{\rm dim}_kM_i}{i^n})$ is an integer in this case (called the multiplicity of $M$). Is $n!({\rm lim sup}_{n\to \infty}\frac{{\rm dim}_kM_i}{i^n})$ an integer in characteristic $p>0$? Does lim$_{n\to \infty}\frac{{\rm dim}_kM_i}{i^n}$ exist in characteristic $p>0$?

\medskip

Since these problems are open only in characteristic $p>0$, it is worth pointing out that Bavula \cite{Ba} has given some striking examples of properties that hold in characteristic 0 but fail in characteristic $p>0$. We briefly mention some of them.

Let a $\mathcal D$-module $M$ be generated by a finite set $z_1\dots,z_s\in M$. Let $M_0$ be the $k$-linear span of $z_1,\dots, z_s$ and let $M_i=\mathcal F_iM_0$. Bavula defines the dimension of $M$ as inf$\{r\in \mathbb R|{\rm dim}_kM_i<i^r\}$ for all sufficiently big $i$. It is not hard to show that this definition is independent of a particular choice of a finite set of generators. In characteristic zero it coincides with the usual definition of the dimension of a finitely generated $\mathcal D$-module.

Bavula shows \cite[9.4]{Ba} that dim$M\geq n$ for every finitely generated $\mathcal D$-module $M$, an analog of the celebrated characteristic zero Bernstein inequality. This inequality is straightforward from \ref{geq}. 

Yet Bavula also shows that there are major differences between characteristic zero and characteristic $p>0$ cases. These are 

(a) in characteristic zero the set of possible values of dim$M$ is all integers between $n$ and $2n$ while in characteristic $p>0$ it is the set of all real numbers between $n$ and $2n$, and 

(b) in characteristic zero a finitely generated $\mathcal D$-module $M$ is holonomic if and only if its dimension is $n$ while in characteristic $p>0$ there exist $M$ such that dim$M=n$ yet $M$ is not holonomic.

3. Perhaps the most interesting open problem is to find a characteristic-free proof of the fact that $\mathcal R_f$ has finite length in the category of $k$-linear $\mathcal D$-modules of the ring $\mathcal R$ of formal power series in a finite number of variables over $k$. Our proof of Theorem \ref{main} suggests that a suitable characteristic-free definition of holonomicity could lead to such a proof.

\end{document}